%

\input amstex
\documentstyle{amsppt}
%
%
\magnification=\magstep1  	
\hsize = 6.5 truein      	
\vsize = 9 truein         	
\TagsAsMath              	
\NoRunningHeads			
%
\topmatter
\title Reverse Mathematics and \\ Recursive Graph Theory
\endtitle
\author  William Gasarch \\
{\rm University of Maryland}\\
Jeffry L. Hirst\\
{\rm Appalachian State University}
\endauthor
\address
Dept. of Computer Science and Institute for Advanced Studies,
University of Maryland,
College Park, MD 20742
\endaddress
\email
gasarch\@cs.umd.edu
\endemail
\address
Dept. of Mathematical Sciences,
Appalachian State University,
Boone, NC  28608
\endaddress
\email
jlh\@math.appstate.edu
\endemail
\date
March 31, 1994
\enddate
\thanks
Gasarch's research was partially supported by NSF grant CCR9020079.
\endthanks
\endkeywords
\endsubjclass
\abstract
We examine a number of results of infinite combinatorics
using the techniques of reverse mathematics.  Our results are
inspired by similar results in recursive combinatorics.
Theorems included
concern colorings of graphs and bounded graphs, Euler paths, and
Hamilton paths.
\endabstract
\endtopmatter
%
\addto\tenpoint{\normalbaselineskip20pt\normalbaselines} 	
\document
\define\rca{\bold {RCA_0}}
\define\wkl{\bold {WKL_0}}
\define\aca{\bold {ACA_0}}
\define\s11{\bold{\Sigma_1^1{-}AC_0}}
\define\atr{\bold {ATR_0}}
\define\p11{\bold {\Pi_1^1{-}CA_0}}

Reverse mathematics provides powerful techniques for analyzing the
logical content of theorems. By contrast, recursive mathematics analyzes
the effective content of theorems.  In many cases, theorems of reverse
mathematics have recursion theoretic corollaries.  Conversely, theorems
and techniques of recursive mathematics can often inspire related results
in reverse mathematics, as demonstrated by the research presented here.
In Section 1, a brief description
of reverse mathematics is given.  Sections 2 and 3 analyze theorems on
graph colorings.  Section 4 considers graphs with Euler paths.  Stronger
axiom systems are introduced in Section 5 and applied to the study
of Hamilton paths in Section 6.

\subhead{1. Reverse mathematics}\endsubhead

In \cite{4},  {\smc Friedman} defined subsystems of second-order arithmetic useful
in determining the proof-theoretic and recursion-theoretic strength of theorems.
The language of second-order arithmetic contains two types of variables,
lower case variables representing elements of $\Bbb N$, the natural numbers,
and upper case variables representing subsets of $\Bbb N$.  Consequently,
a model for a subsystem consists of a number universe and a collection of
subsets of the number universe.

$\rca$ is the weak base system used in reverse mathematics.  It consists of
the axioms of first order Peano arithmetic with induction restricted to
$\Sigma_1^0$ formulas, and the recursive comprehension axiom, which states
that any set definable by both a $\Sigma_1^0$ formula and a $\Pi_1^0$ formula
exists.  $\rca$ suffices to prove fundamental facts about pairing functions, finite
sequences, and other tools used to encode theorems as statements of second-order
arithmetic.  In this paper, much of the coding has been suppressed.  Details on encoding
techniques can be found in \cite{16}.

Stronger axiom systems can be constructed by adding additional set existence axioms
to $\rca$.  For example, the subsystem $\wkl$ consists of the axioms of $\rca$ together
with a weak version of K\"onig's Lemma asserting that every infinite $0$--$1$ tree contains an
infinite path. $\wkl$ is strictly stronger than $\rca$.  Often it is possible to show that a theorem is
equivalent to a set comprehension axiom over the weak base system $\rca$.  Results of this sort, called reverse mathematics, leave no doubt as to what set existence axioms are necessary
in a proof.  The following theorem of {\smc Simpson} \cite{14} illustrates this process, and is used
in later sections.
The notation ($\rca$) in the proclamation of a theorem or definition signifies that the
theorem can be proved in $\rca$, or that the definition can be expressed in the language
of $\rca$ using coding techniques.

\proclaim{Theorem 1 ($\rca$)}
The following are equivalent:
\roster
\item  $\wkl$.
\item  If $f : \Bbb N \rightarrow \Bbb N$ and $g: \Bbb N \rightarrow \Bbb N$ are
injections such that for
all $j,k \in \Bbb N$, $f(j) \neq g(k)$, then there is a set $X$ that separates
the ranges of $f$ and $g$, formally,
$$\forall j \forall n((f(j)=n \rightarrow n \in X) \land (g(j)=n \rightarrow n \notin X )).$$
\endroster
\endproclaim

Adopting a model theoretic viewpoint can clarify the content of Theorem 1.
In part, the theorem asserts that if $f$ and $g$ are injections (encoded) in
a model of $\wkl$, then a separating set for $f$ and $g$ is also (encoded)
in the model.   In some sense, this implicitly restricts the choices of $f$ and $g$.

The axiom system $\aca$ consists of $\rca$ together with the arithmetical comprehension
scheme.  This scheme asserts that any set definable by a formula containing no set quantifiers
exists.  $\aca$ is strictly stronger than $\wkl$.  A proof of the following characterization
of $\aca$ can be found in {\smc Simpson} \cite{14}.

\proclaim{Theorem 2 ($\rca$)}
The following are equivalent:
\roster
\item  $\aca$.
\item  If $f: \Bbb N \rightarrow \Bbb N$ is an injection, then the range of $f$ exists.
\endroster
\endproclaim

Additional axiom systems are briefly described in section 5.  For more detailed information
on subsystems of second-order arithmetic and reverse mathematics, see \cite{15} or \cite{16}.

\subhead{2. Graph Colorings}\endsubhead

In this section we will consider theorems on node colorings of countable graphs.
A (countable) graph $G$ consists of a set of vertices $V \subseteq \Bbb N$ and
a set of edges $E \subseteq [\Bbb N ]^2$.
We will abuse notation by denoting an edge by $(x,y)$ rather than
$\{ x,y \}$.
For $k \in \Bbb N$, we say
that $\chi : V \rightarrow k$ is a $k$-coloring of $G$ if
$\chi$ always assigns different colors to neighboring vertices.
That is, $\chi$ is a $k$-coloring if $\chi : V \rightarrow k$ and
$(x,y) \in E$ implies $\chi (x) \neq \chi (y)$.  If $G$ has a $k$-coloring,
we say that $G$ is $k$-chromatic.
Using an appropriate axiom system, it is possible to prove that a graph
is $k$-chromatic if it satisfies the following local condition.

\definition{Definition 3 ($\rca$)}
A graph $G$ is
{\it locally $k$-chromatic}
if every finite subgraph of $G$ is $k$-chromatic.
\enddefinition

The following theorem is the simplest result concerning graph colorings.
To prove that \therosteritem1 implies \therosteritem2, a tree is constructed
in which every infinite path encodes a $k$-coloring.  The proof of the reversal
uses a graph whose $k$-colorings encode separating sets for a pair of injections.
Theorem 1 is then applied to finish the proof.  For a detailed proof, see Theorem 3.4
in \cite{9}.

\proclaim{Theorem 4 ($\rca$)}
For every $k \ge 2$, the following are equivalent:
\roster
\item  $\wkl$.
\item  If $G$ is locally $k$-chromatic, then
$G$ is $k$-chromatic.
\endroster
\endproclaim

>From Theorem 4, we can deduce two recursion theoretic results due to {\smc Bean} \cite{2}.
The first result can be proved directly by imitating the construction used in the proof
of the reversal of Theorem 4, using a pair of recursive functions with no recursive
separating set.  We will provide an alternative model theoretic argument based on
the following observation.  
By Theorem 1
and the existence of a pair of recursive functions with no recursive separating set,
every $\omega$-model of $\wkl$ must contain a non-recursive set.

\proclaim{Corollary 5 ({\smc Bean} \cite{2}) }
For every $k \ge 2$, there is a recursive $k$-chromatic
graph which has no recursive $k$-coloring.
\endproclaim
\demo{Proof}
Suppose, by way of contradiction, that for some $k \in \omega$
every recursive $k$ chromatic graph has a recursive $k$-coloring.
Then, by Theorem 4,
$\omega$ together with the recursive sets is a model
of $\wkl$, contradicting the fact that every $\omega$-model of
$\wkl$ contains a non-recursive set.  \qed
\enddemo

Our model theoretic proof of the
next recursion theoretic corollary relies on the fact that the set universe
of any $\omega$-model of  $\wkl$ is a Scott system \cite{12}.
Such a model will include the recursive sets, and additional sets which can
be bounded in complexity.  By the Shoenfield-Kreisel low basis theorem
\cite{13} there is an $\omega$-model of $\wkl$ such that for each set
$X$ in the model, if $a$ is the Turing degree of $X$, then
$a^\prime \le 0^\prime$.  That is, every set in such a model of $\wkl$ is of
low degree.

\proclaim{Corollary 6 ({\smc Bean} \cite{2})}
For every $k \ge 2$, every recursive $k$-chromatic graph
has a $k$-coloring of low degree.
\endproclaim
\demo{Proof}
Let $M$ be an $\omega$-model of $\wkl$ in which every set is of low degree.
Let $G$ be a recursive $k$-chromatic graph.  Then $G$ is (encoded) in $M$,
and by Theorem 4, a $k$-coloring of $G$ is also (encoded) in $M$.  Thus $G$
has a $k$-coloring of low degree.  \qed
\enddemo

The number of colors allowed in a coloring of a locally $k$-chromatic graph
can be substantially increased without weakening the logical strength of the
resulting theorem.  This contrast sharply with the situation for bounded graphs
which is discussed in the next section.

\proclaim{Theorem 7 ($\rca$)}
For each $k \ge 2$, the following are equivalent:
\roster
\item $\wkl$.
\item If $G$ is locally $k$-chromatic, then $G$ is
$(2k-1)$-chromatic.
\endroster
\endproclaim
\demo{Proof}
Since $\rca$ proves that every $k$-chromatic graph is
$(2k-1)$-chromatic, \therosteritem1 implies
\therosteritem2 follows immediately
from
Theorem 4.

We will now prove that \therosteritem2
implies \therosteritem1 when $k=2$, and
then indicate how the argument can be generalized to any
$k \in \Bbb N$.
By Theorem 1, $\wkl$ can be proved by showing that the ranges
of an arbitrary pair of disjoint injections can
be separated.  Let
$f : \Bbb N \rightarrow \Bbb N$ and
$g : \Bbb N \rightarrow \Bbb N$ be
injections such that for all $m, n \in \Bbb N$,
$f(n) \neq g(m)$.  We will construct a $2$-chromatic
graph with the property that any $3$-coloring of
$G$ encodes a set $S$ such that $y \in Range (f)$
implies $y \in S$, and $y \in S$ implies
$y \notin Range (g)$.

The graph $G$ contains an infinite complete bipartite
subgraph consisting of upper vertices
$\{ b^u_n : n \in \Bbb N \}$, lower vertices
$\{ b^l_n : n \in \Bbb N \}$, and connecting
edges $\{ ( b^u_n , b^l_m ) : n,m \in \Bbb N \}$.
Also, $G$ contains an infinite collection of pairs of vertices,
denoted by $n^u$ and $n^l$ for $n \in \Bbb N$.
Each such pair is connected, so the edges
$\{ ( n^u , n^l ) : n \in \Bbb N \}$ are included in $G$.
Additional connections depend on the injections $f$ and $g$.
If $f(i)=n$, add the edges $(b^u_m, n^l)$ and
$(b^l_m, n^u)$ for all $m \ge i$.  If
$g(i)= n$, add the edges $(b^u_m, n^u)$ and
$(b^l_m , n^l)$ for all $m \ge i$.
Naively, if $n$ is in the range of $f$ or $g$, then the
pair $(n^u, n^l)$ is connected to the complete bipartite
subgraph.  If $n$ is in the range of $G$, the pair is
``flipped'' before it is connected.
The reader can verify that $G$ is $\Delta_1^0$ definable in
$f$ and $g$, and thus exists by the recursive comprehension axiom.
Every finite subgraph of $G$ is clearly bipartite,
so $G$ is locally $2$-chromatic.  Thus, by (2), $G$ has
a 3-coloring; denote it by $\chi : G \rightarrow 3$.

If $\chi$ is a $2$-coloring, we can define the separating set, $S$, by
$$
S = \{ y \in \Bbb N : \chi (y^u) = \chi (b_0^u) \lor
\chi (y^l) = \chi ( b_0^l) \}.
$$
When $\chi$ uses all $3$ colors, we must modify the construction of $S$.
In particular, we must find a $j \in \Bbb N$ such that
\roster
\item"({\it a})"
$\forall y ( \exists n ( n\ge j \land f(n)=y) \rightarrow
(\chi (y^u) = \chi (b^u_j) \lor \chi (y^l)= \chi (b^l_j)))$, and
\item"({\it b})"
$\forall y ( \exists n ( n\ge j \land g(n) = y ) \rightarrow
(\chi ( y^l ) \neq \chi ( b^l_j ) \land
\chi (y ^u ) \neq \chi ( b_j^u )))$.
\endroster
Suppose, by way of contradiction, that no such $j$ exists.
The for some $m$ and $y$, either
$f(m) = y \land \chi (y^u ) \neq \chi(b^u_0) \land \chi (y^l) \neq \chi(b^l_0)$
or
$g(m) = y \land (\chi(y^l)=\chi(b^l_0)\lor\chi(y^u)=\chi(b^u_0))$.
If $f(m) = y$, since $\chi$ is a $3$-coloring,
either $\chi(y^u)=\chi(b^l_0)$ or
$\chi(y^l) = \chi(b_0^u)$.
By the construction of $G$, for every $n > m$,
$\chi(b^u_m)=\chi(b^u_n)$ and
$\chi(b^l_m)=\chi(b^l_n)$.
Similarly, the case $g(m)=y$ also yields a point
beyond which the complete bipartite subgraph of $G$ is $2$-colored.
By the negation of ({\it a}) and ({\it b}), there
is an $m^\prime > m$ and a $z \in \Bbb N$
such that either
$f(m^\prime) = z \land \chi (z^u ) \neq \chi(b^u_m)
\land \chi (z^l)\neq\chi(b^l_m)$
or
$g(m^\prime) = z \land (\chi(z^l)=\chi(b^l_m)\lor\chi(z^u)=\chi(b^u_m))$.
If $f(m^\prime )=z$, then since $\chi$ is a $3$-coloring,
either $\chi(z^u)=\chi(b^l_m)$ or
$\chi(z^l) = \chi(b^u_m)$.  Since $m^\prime > m$,
$\chi(b^l_m)=\chi(b^l_{m^\prime})$ and
$\chi(b^u_m)=\chi(b^u_{m^\prime})$, so either
$\chi(z^l)=\chi(b^u_{m^\prime})$ or
$\chi(z^u)=\chi(b^l_{m^\prime})$.
But $(z^l,b^u_{m^\prime})$ and
$(z^u,b^l_{m^\prime})$ are edges of $G$, so $\chi$ is not a
$3$-coloring.  Assuming $g(m^\prime )=z$ yields a similar
contradiction.  Thus, a $j$ satisfying ({\it a}) and ({\it b})
exists.

Given an integer $j$ satisfying ({\it a}) and ({\it b}),
the separating set $S$ may be defined as the union of
$\{ y\in \Bbb N : \exists n<j \, f(n) = y \}$ and
$$
\{
y \in \Bbb N :
( \forall n<j \, g(n) \neq y)
\land
(\chi (y^u) = \chi (b^u_j) \lor \chi (y^l) = \chi(b^l_j))
\}
$$
$S$ is $\Delta_1^0$ definable in $\chi$ and $j$,
so the recursive comprehension axiom assures the existence of
$S$.  If $f(n) = y$ and $n<j$, then $y \in S$.
If $f(n) = y$ and $n \ge j$, then by ({\it a}) and
the fact that $f$ and $g$ have disjoint ranges, $y \in S$.
Thus $Range (f) \subseteq S$.  If $g(n)=y$, and $n<j$, then since
the ranges of $f$ and $g$ are disjoint we have $y \notin S$.
If $g(n) = y$ and $n \ge j$, by ({\it b}) $y \notin S$.
Thus $S$ is the desired separating set.  This completes the
proof for $k=2$.

For $k>2$, the preceding proof requires the following
modifications.  Replace the complete bipartite subgraph
of $G$ by a complete $k$-partite subgraph with vertices
$\{b^p_m: p<k \land m\in \Bbb N \}$.  Each pair
$( n^u , n^l )$ is replaced by a complete graph on
the vertices $\{ n^p : p<k \}$.
If $f(i) = n$, add the edges
$(b^p_m, n^{p^\prime})$ for all $m \ge i$ and
all $p \neq p^\prime$ less than $k$.
If $g(i) = n$, twist the subgraph before attaching it.  That is,
add the edges
$( b^p_m , n^ {p^\prime} )$ for all $m \ge i$ and all
$p$ and $p^\prime$ less than $k$ such that
$p \not\equiv p^\prime +1$ (mod $k$).
The argument locating the integer $j$ is similar,
except that $m$ and $m^\prime$ must be replaced
by a sequence $m_1 , \dots , m_k$.  Beyond the
point $m_{k-1}$, the complete $k$-partite subgraph
of $G$ is $k$-colored by $\chi$.  The definition of
$S$ is very similar, except that a bounded quantifier
should be used to avoid the $k$-fold conjunction.  \qed
\enddemo

The following recursion theoretic consequence of Theorem 7 is a
special case of a result due to {\smc Bean}.

\proclaim{Corollary 8 ({\smc Bean} \cite{2})}
For every $k \ge 2$, there is a recursive graph $G$ which has
no recursive $(2k{-}1)$-coloring.
\endproclaim

\demo{Proof}
Imitate the reversal of Theorem 7, using disjoint recursive injections
with no recursive separating set.  \qed
\enddemo

{\smc Bean} \cite{2} showed that Corollary 8 holds with $2k-1$ replaced
by any value larger than $k$.  In light of this, the following conjecture
seems reasonable.  Unfortunately, even the case where $k=2$ and $m=4$
remains open.

\proclaim{Conjecture 9 ($\rca$)}
For each $k \ge 2$ and each $m \ge k$ the following are equivalent:
\roster
\item $\wkl$.
\item If $G$ is locally $k$-chromatic, then $G$ is
$m$-chromatic.
\endroster
\endproclaim

\remark{Remark}
Note that \therosteritem1 implies \therosteritem2 follows from Theorem 7.
Also, the full reversal is an easy corollary of the reversal for
$k=2$ and arbitrary $m$.  To see this, note that
if $G$ is the graph used to
prove the reversal for $k$ and $m$, the graph resulting from adding
one vertex to $G$ and attaching it to every existing vertex will provide
a proof of the reversal for $k+1$ and $m+1$.
\endremark

\subhead{3. Bounded graphs and sequences of graphs}\endsubhead

As noted above, a locally $k$-chromatic recursive graph may not
have a recursive coloring, regardless of the number of colors used.
By contrast, highly recursive graphs always have recursive colorings.
A proof theoretic
analog of a highly recursive graph is a bounded graph.

\definition{Definition 10 ($\rca$)}
A graph $G = \langle V,E \rangle$ is
{\it bounded}
if there is a function $h : V \rightarrow \Bbb N$ such that
for all $x, y \in V$, $(x,y) \in E$ implies $h(x) \ge y$.
\enddefinition

Using Definition 10, we can state a proof theoretic version of a
theorem on highly recursive graphs proved by {\smc Schmerl} \cite{11} and
independently rediscovered by {\smc Carstens} and {\smc Pappinghaus} \cite{3}.

\proclaim{Theorem 11 ($\rca$)}
For $k \in \Bbb N$, if $G$ is a bounded locally $k$-chromatic graph,
then $G$ is $({2k-1})$-chromatic.
\endproclaim

\demo{Proof}
The proof of Theorem 1 of {\smc Schmerl} \cite{11} can be carried out in $\rca$.  \qed
\enddemo

\proclaim{Corollary 12 ({\smc Schmerl} \cite{11})} 
For every $k$, every highly recursive $k$-chromatic graph
has a recursive $({2k-1})$-coloring.
\endproclaim

If the number of colors allowed is less than $2k-1$, a $k$-chromatic
highly recursive graph
may not have a recursive coloring.

\proclaim{Theorem 13 ($\rca$)}
For every $k \ge 2$, the following are equivalent:
\roster
\item $\wkl$.
\item
If $G$ is a bounded locally $k$-chromatic graph, then $G$ is
$({2k-2})$-chromatic.
\endroster
\endproclaim

\demo{Proof}
Since $\rca$ proves that every $k$-chromatic graph is
$({2k-2})$-chromatic, \therosteritem1 implies \therosteritem2 follows immediately
from
Theorem 4.  For the case $k=2$, the statement that
\therosteritem2 implies \therosteritem1 is included in Theorem 3.4 of \cite{9}.
Alternately, this case could be proved by formalizing the proof of
Theorem 2 of {\smc Schmerl} \cite{11}.  We will adopt this approach for the case $k>2$.

Let $k>2$, and assume that \therosteritem2 holds.
As in
Theorem 7, we will prove $\wkl$ by finding a separating set for the
ranges of two disjoint injections, $f$ and $g$.  The separating set
must be encoded in a $({2k-2})$-coloring of a bounded locally $k$-chromatic
graph $G$.

$G$ will be constructed from subgraphs called {\it blocks}.
A block $B$ consists of $k^2$ vertices
$\{ v_{ij} : i<k \land j<k \}$ connected by the edges
$(v_{ij}, v_{rs})$ for $i \neq r$ and $j \neq s$.
A block can be viewed as a $k \times k$ matrix where each
vertex is connected to all the elements of its associated
cofactor matrix.  We can {\it link} two blocks $B$ and $B^\prime$ by
adding the edges $( v_{ij} , v_{rs}^\prime )$ for $i \neq r$ and
$j \neq s$.

Given a coloring $\chi$ of a block $B$, we say that $B$
{\it has a colorful row} if for some $i$, whenever
$j \neq r$, $\chi ( v_{ij} ) \neq \chi ( v_{ir} ) $.
Similarly, $B$ has a colorful column if all the elements in some column
have distinct colors.  $\rca$ proves that if $\chi$ is
a $({2k-2})$-coloring of a block $B$, then $B$ has either
a colorful row or a colorful column, but not both.
(To prove this, formalize the proof of Lemma 2.1 in \cite{11}
or Lemma 5.25 in \cite{6}.)  Furthermore, $\rca$ proves
that if $\chi$ is a $({2k-2})$-coloring of two linked blocks
$B$ and $B^\prime$, then $B$ has a colorful row (column)
if and only if $B^\prime$ has a colorful column (row).
(To prove this, formalize the proof of Lemma 2.2 in \cite{11}
or Lemma 5.27 in \cite{6}).

The graph $G$ is constructed from two sets of blocks,
$\{B_j : j \in \Bbb N \}$ and $\{B_{ij} : i,j \in \Bbb N \}$.
For each $i$ and $j$, $B_j$ is linked to $B_{(j+1)}$, and
$B_{ij}$ is linked to $B_{i(j+1)}$.  Additional links depend
on the injections $f$ and $g$.  If $f(m) = n$, link
$B_{n(2m)}$ to $B_{(2m)}$.  If $g(m)=n$, link $B_{n(2m)}$ to
$B_{(2m+1)}$.  The reader may verify that $G$ is $\Delta_1^0$ definable
in $f$ and $g$, bounded, and locally $k$-chromatic.  Applying
\therosteritem2, $G$ has a $({2k-2})$-coloring, $\chi$.
By the recursive comprehension axiom, the set
$$S = \{ n: B_{n0}~\text{has a colorful column}\}$$
exists.  We will show that $S$ is the desired separating set.
Suppose first that $\chi$ induces a colorful row in $B_0$.
If $f(m)=n$, since $B_0$ and $B_{n0}$ are connected by a sequence of
linked blocks of even length, $B_{n0}$ has a colorful column, and $n \in S$.
Also, if $g(m) = n$, $B_0$ and $B_{n0}$ are linked by an odd length
sequence, so $n \notin S$.  Thus, $Range (f) \subseteq S$ and
$Range (g) \cap S  = \emptyset$.  Similarly, if $\chi$ induces a colorful
column in $B_0$, then $S$ is a separating
set containing $Range (g)$.  \qed
\enddemo

Imitating the reversal of Theorem 13 using disjoint recursive injections with
recursively inseparable ranges yields {\smc Schmerl's} proof of the following result.

\proclaim{Corollary 14 ({\smc Schmerl} \cite{11})}
For each $k \ge 2$ there is a highly recursive $k$-chromatic graph
which has no recursive ${2k-2}$ coloring.
\endproclaim

We will close this section with a theorem concerning {\sl sequences \/} of graphs
and its recursion theoretic corollary.
We say that a graph $G$ is {\it colorable} if
there exists an integer $k$ such that $G$ is $k$-chromatic.

\proclaim{Theorem 15 ($\rca$)}
The following are equivalent:
\roster
\item
$\aca$.
\item
Given a countable sequence of graphs,
$\langle G_i : i \in \Bbb N \rangle$,
there is a function $f : \Bbb N \rightarrow 2$
such that $f(i) = 1$ if $G_i$ is colorable and
$f(i)=0$ otherwise.
\endroster
\endproclaim

\demo{Proof}
To prove that \therosteritem1 implies \therosteritem2, assume
$\aca$ and let $\langle G_i : i \in \Bbb N \rangle$ be
a sequence of graphs.  Define $f : \Bbb N \rightarrow \Bbb N$
by setting $f(i)=1$ if there exists a $k \in \Bbb N$ such that
$G_i$ is locally $k$-chromatic, and setting $f(i)=0$ otherwise.
Since ``$G_i$ is locally $k$-chromatic'' is an arithmetical
sentence with parameter $G_i$,  $f$ exists by the arithmetical
comprehension axiom.  Since $\aca$ implies $\wkl$, we may
apply Theorem 4 to show that $f(i) = 1$ if and only if $G_i$ is
colorable.

To prove the converse, assume $\rca$ and \therosteritem2.  By
Theorem 2, to
prove $\aca$ it suffices to show that for every injection $g$,
$Range (g)$ exists.
Define the sequence of graphs $\langle G_i : i \in \Bbb N \rangle$
as follows.  Let $\{ v_j : j \in \Bbb N \}$ be the vertices
of $G_i$.  If $j<k$ and
$\forall m \le k ( g(m) \neq i)$, add the edge $(v_j , v_k )$ to
$G_i$.  $\rca$ can prove that $\langle G_i : i \in \Bbb N \rangle$
exists, and $G_i$ is colorable if and
only if $i \in Range (g) $.
Thus, the function $f$ supplied by \therosteritem2 is the characteristic
function for $Range (g)$.  By the recursive comprehension axiom,
$Range (g)$ exists.  \qed
\enddemo

\proclaim{Corollary 16}
There is a recursive sequence of recursive graphs
$\langle G_i : i \in \Bbb N \rangle$ such that
$0^\prime$ is recursive in
$\{ i \in \Bbb N : G_i ~\text{is colorable}\}$.
\endproclaim

\demo{Proof}
In the proof of the reversal for Theorem 15,
let $g$ be a recursive function such that $0^\prime$ is recursive
in $Range (g)$.  The sequence of graphs constructed in
the proof has the desired properties.  \qed
\enddemo

\subhead{4. Euler paths}\endsubhead

Now, we will turn to the study of Euler paths.
A {\it path} in a graph $G$ is a sequence of vertices
$v_0 , v_1 , v_2 , \dotsc$ such that for every $i \in \Bbb N$,
$( v_i , v_{i+1} )$ is an edge of $G$.
A path is called an {\it Euler path} if it uses every edge of
$G$ exactly once.

The following terminology is useful in determining when a graph
has an Euler path.
A graph $G = \langle V , E \rangle$ is {\it locally finite}
if for each vertex $V$, the set $\{ u \in V : (v,u) \in E \}$ is finite.
If $H$ is a subgraph of $G$, $G - H$ denotes the graph obtained by
deleting the edges of $H$ from $G$.  Using this terminology,
we can describe a condition which, from a naive viewpoint, is
sufficient for the existence of an Euler path.

\definition{Definition 17 ($\rca$)}
A graph $G$ is
{\it pre-Eulerian}
if it is
\roster 
\item
connected,
\item
has at most one vertex of odd degree,
\item
if it has no vertices of odd degree, then it has at least one vertex of infinite degree, and
\item
if $H$ is any finite subgraph of $G$ then $G - H$ has exactly one
infinite connected component.
\endroster
\enddefinition

Note that the formula ``$G$ is pre-Eulerian'' is arithmetical in the
set parameter $G$.  $\rca$ suffices to prove that every graph with
an Euler path is pre-Eulerian.  However, $\rca$ can only prove that
bounded pre-Eulerian graphs have Euler paths.
(Bounded graphs are defined in Section 3.)
This result is just a formalization of
{\smc Bean's} \cite{1} proof that every highly recursive
pre-Eulerian graph has a recursive Euler path.

\proclaim{Theorem 18 ($\rca$)}
If $G$ is a bounded pre-Eulerian graph, then $G$ has an Euler path.
\endproclaim

\demo{Proof}
The proof of this theorem is just a straightforward formalization
of Theorem 2 of {\smc Bean} \cite {1}.  The formalization requires verification
that Euler's Theorem for finite graphs (see \cite{10}) can be proved
using only $\rca$.  \qed
\enddemo

If $G$ is not bounded, additional axiomatic strength is required to prove
the existence of an Euler path.

\proclaim{Theorem 19 ($\rca$)}
The following are equivalent:
\roster
\item 
$\aca$
\item
If $G$ is a pre-Eulerian graph, then $G$ has an Euler path.
\item
If $G$ is a locally finite pre-Eulerian graph, then $G$ has an
Euler path.
\endroster
\endproclaim

\demo{Proof}
To prove that \therosteritem1 implies \therosteritem2,
assume $\aca$ and let $G$ be a pre-Eulerian graph.
Let $\langle E_i : i \in \Bbb N \rangle$ be an enumeration of
the edges of $G$.  Let $v_0$ be the vertex of $G$ of odd degree,
or a vertex of infinite degree if no odd vertex exists.
Imitating the proof of Theorem 3.2.1 of Ore \cite{10},
there is a finite path $P$ containing the edge $E_0$ such that
\roster
\item"$\bullet$"
$P$ starts at $v_0$,
\item"$\bullet$"
$G-P$ is connected, and
\item"$\bullet$"
$P$ ends at the odd vertex of $G-P$, or at an infinite vertex
of $G-P$ if no odd vertex exists.
\endroster
Furthermore, since the finite paths of $G$ can be encoded
by integers, we can pick the unique path $P_0$ satisfying
the conditions above and having the least code.  Similarly,
any path $P_i$ satisfying the three conditions can be extended
to a unique path $P_{i+1}$ which contains the edge $E_{i+1}$,
satisfies the three conditions, and has the least code among all
paths with these properties.  Note the $P_{i+1}$ extends $P_i$
by including $P_i$ as an initial segment.  The reader
may verify that the sequence of paths $\langle P_i : i \in \Bbb N \rangle$
is arithmetically definable in $G$, and so exists by arithmetical comprehension.
Let $v_{i}$ denote the $i^{th}$ vertex of $P_i$.  Then the
sequence $\langle v_{i} : i \in \Bbb N \rangle$ exists by recursive
comprehension and includes each $P_i$ as an initial segment.
Consequently, $\langle v_{i}: i \in \Bbb N \rangle$ defines
an Euler path through $G$.

Since \therosteritem3 is a special case of \therosteritem2,
showing that \therosteritem3 implies \therosteritem1 will complete
the proof of the theorem.
Assume $\rca$ and fix an injection $f: \Bbb N \rightarrow \Bbb N$.
We will construct a locally finite pre-Eulerian
graph $G$ such that every Euler path
through $G$ encodes $Range (f)$.  Define the vertices of $G$ by
$$
V = \{ a_n , b_n , c_n : n \in \Bbb N \}.
$$
For each $n$, include the edges $( a_n , a_{n+1} )$ and
$( b_n , c_n )$ in $G$.  Additionally, for each
$i$ and $n$, if $f(i)=n$ then include the edges
$( a_n , b_i )$ and $ ( c_i , a_n )$
in $G$.  $\rca$ suffices to prove that
$G$ exists, and is both locally finite and pre-Eulerian.
By \therosteritem3, $G$ has an Euler path.
Note that $n \in Range(f)$ if and only if
the first occurrence of $a_n$ in the Euler path is not
followed immediately by $a_{n+1}$.  By the recursive
comprehension axiom, $Range(f)$ exists.  Since $f$ was
an arbitrary injection, by Theorem 2 this suffices to prove
$\aca$.  \qed
\enddemo

\proclaim{Corollary 20}
There is a recursive pre-Eulerian graph $G$ such that
$0^\prime$ is recursive in every Euler path through $G$.
\endproclaim

\demo{Proof}
Let $f$ be a recursive function such that $0^\prime$ is
recursive in $Range (f)$.  Construct the graph $G$ as
in the proof of the reversal in Theorem 19.
Then $G$ is recursive, and $Range (f)$ is recursive
in every Euler path through $G$.  \qed
\enddemo

$\aca$ also suffices to address the problem of determining which
elements of a sequence of graphs have Euler paths.  This contrasts
sharply with the situation for Hamilton paths, as described in Theorem 30.

\proclaim{Theorem 21 ($\rca$)}
The following are equivalent:
\roster
\item 
$\aca$
\item
Given a countable sequence of graphs,
$\langle G_i : i \in \Bbb N \rangle$,
there is a set $Z \subseteq \Bbb N$ 
such that $i \in Z$ if and only if
$G_i$ has an Euler path.
\endroster
\endproclaim

\demo{Proof}
First assume $\aca$ and let
$\langle G_i : i \in \Bbb N \rangle$
be a sequence of graphs.  Define the set $Z$ by 
$Z = \{ i \in \Bbb N : G_i ~ \text{is pre-Eulerian} \}$.
Note that $Z$ is arithmetically definable in
$\langle G_i : i \in \Bbb N \rangle$, so $\aca$ proves the existence of $Z$.
Since $\rca$ proves that every graph with an Euler path is pre-Eulerian,
and Theorem 19 proves that every pre-Eulerian graph has an Euler path, 
$i \in Z$ if and only if $G_i$ has an Euler path.

To prove the converse, assume $\rca$ and \therosteritem2.
By Theorem 2,
it suffices to prove that $Range(f)$ exists for an arbitrary
injection $f: \Bbb N \rightarrow \Bbb N$.
Let
$v_0 , v_1 , v_2 , \dotsc$ denote
the vertices of $G_i$.  For each $n$,
if $f(n) \neq i$, add the
edge $( v_n , v_{n+1} )$ to $G_i$.
By the recursive comprehension axiom, the
sequence of graphs $\langle G_i : i \in \Bbb N \rangle$
exists.  Let $Z$ be as in \therosteritem2.
Then
$Range (f) = \{ i \in \Bbb N : i \notin Z \}$,
so $Range (f)$ exists by the recursive comprehension axiom.
Note that this proof actually shows that Theorem 21 holds with
\therosteritem2 restricted to sequences of bounded graphs.
\qed
\enddemo

Theorem 21 can be used to establish rough upper and lower
bounds for the complexity of the problem of determining which
graphs in a sequence have Euler paths.

\proclaim{Corollary 22 ({\smc Beigel} {\rm and} {\smc Gasarch} \cite{6})}
If $\langle G_i : i \in \Bbb N \rangle$ is an arithmetical
sequence of graphs, then the set
$\{ i \in \Bbb N : G_i ~\text{has an Euler path} \}$ is arithmetical.
\endproclaim

\demo{Proof}
$\omega$ together with the arithmetical sets is a model of
$\aca$, and thus models \therosteritem2 of Theorem 21.  \qed
\enddemo

\proclaim{Corollary 23 ({\smc Beigel} {\rm and} {\smc Gasarch} \cite{6})}
There is a recursive sequence of recursive graphs,
$\langle G_i : i \in \Bbb N \rangle$
such that $0^\prime$ is recursive in the set
$\{ i \in \Bbb N : G_i ~\text{has an Euler path} \}$.
\endproclaim

\demo{Proof}
Let $f$ be a recursive function such that $0^\prime$ is recursive
in $Range (f)$.  The sequence of graphs constructed from $f$ as in
the proof of the reversal in Theorem 21 has the desired property.  \qed
\enddemo

\remark{Remark}
A {\it two-way} or {\it endless} Euler path is a bijection between the
integers (both positive and negative) and the set of edges of $G$ such
that  each edge shares one vertex with its predecessor and its other vertex
with its successor.  Theorems 18, 19, and 21 can be modified to address
the existence of two-way Euler paths.
\endremark

\subhead{5. Stronger axiom systems}\endsubhead

The discussion of Hamilton paths in the next section uses three axiom systems
which are each strictly stronger than $\aca$.  These axiom systems, in strictly
increasing order of strength, are
$\s11$, $\atr$, and $\p11$.

The subsystem $\Sigma_1^1$-axiom of choice, denoted by
$\s11$, consists of $\rca$ together with the
comprehension scheme
$$
(\forall k ( \exists X \, \Psi (k,X))) \rightarrow (\exists Y ( \forall k \, \Psi (k,(Y)_k)))
$$
where $\Psi$ is any $\Sigma _1^1$ formula and
$(Y)_k = \{ i : (i,k) \in Y \}$.

The subsystem $\atr$ consists of $\aca$ and
an existence axiom for sets constructed
by applying a form of arithmetical transfinite recursion.
We will need the following notation.  Let $Seq$ denote the set
of (codes for) finite sequences of elements of $\Bbb N$. Given
$T \subseteq Seq$, we say that $T$ is a {\it tree} if
whenever $\tau \in T$ and $\sigma$ is an initial segment of $\tau$,
$\sigma \in T$.  In this way, we can encode infinitely splitting trees
as subsets of $Seq$, which can in turn be encoded as subsets of $\Bbb N$.
The following result (which is Theorem 5.2 of \cite{16})
gives two combinatorial characterizations of $\atr$.

\proclaim{Theorem 24 ($\rca$)}
The following are equivalent
\roster
\item
$\atr$.
\item
The schema
$$
( \forall i )( \exists ~\text{at most one} ~X ) \Psi (i,X) \rightarrow
( \exists Z)( \forall i )( i \in Z \leftrightarrow (\exists X ) \Psi (i,X )) ,
$$
where
$\Psi (i, X)$ is any arithmetical formula in which $Z$ does not occur.
\item
For any sequence of trees
$\langle T_i : i \in  \Bbb N \rangle$, if for each $i \in \Bbb N$,
$T_i$ has at most one path, then
$\exists Z \forall i (i \in Z \leftrightarrow T_i ~\text{has a path})$.
\endroster
\endproclaim

The system $\p11$ consists of $\rca$ plus a comprehension axiom
asserting that the set $\{ n \in \Bbb N : \Psi (n) \}$ exists
for any $\Pi_1^1$ formula $\Psi$.
$\p11$ is strictly
weaker than $\bold {\Pi^1_\infty{-}CA_0}$ (full second-order
arithmetic.) 

\subhead{6. Hamilton paths}\endsubhead

Now we will consider theorems on the existence of
Hamilton paths.
A path through a graph $G$ is called a (one way)
{\it Hamilton path} if it uses every
vertex of $G$ exactly once.
There is no known analog of the characterization
``pre-Eulerian'' for graphs containing Hamilton paths.
Consequently, all the results of this section concern sequences
of graphs.

Each reversal in this section will rely on the construction
of a sequence of graphs from a sequence of trees, as
in the following lemma.

\proclaim{Lemma 25 ($\rca$)}
Given a sequence of trees
$\langle T_i : i \in \Bbb N \rangle$,
there is a sequence of graphs
$\langle G_i : i \in \Bbb N \rangle$
such that
\roster
\item 
for each $i \in \Bbb N$, $T_i$ has a (unique) path if and
only if $G_i$ has a (unique) Hamiltonian path, and
\item
if there is a sequence
$\langle P_i : i \in \Bbb N \rangle$
such that $P_i$ is a Hamiltonian path through
$G_i$ for each $i \in \Bbb N$, then there is
a sequence
$\langle P_i^{\prime} : i \in \Bbb N \rangle$
such that $ P_i^{\prime}$ is a path through
$T_i$ for each $i \in \Bbb N$.
\endroster
\endproclaim

\demo{Proof}
For each $T_i$, use the graph constructed in the proof of
Theorem 1 of {\smc Harel} \cite {5}. \qed
\enddemo

The next three theorems analyze the following tasks:
\roster
\item
finding Hamilton paths through graphs known to have such paths,
\item
determining whether graphs that have at most one Hamilton path have such a path, and
\item
determining whether arbitrary graphs have Hamilton paths.
\endroster
Using proof theoretic strength as a measure of difficulty, we shall see that
these tasks are strictly increasing in order of difficulty.

\proclaim{Theorem 26 ($\rca$)}
The following are equivalent:
\roster
\item 
$\s11$.
\item
If $\langle G_i : i \in \Bbb N \rangle$
is a sequence of graphs such that each $G_i$ has
a Hamilton path, then there is a sequence
$\langle P_i : i \in \Bbb N \rangle$ such
that for each $i$, $P_i$ is a Hamilton path
through $G_i$.
\endroster
\endproclaim

\demo{Proof}
To prove that \therosteritem1 implies \therosteritem2,
assume $\s11$ and let
$\langle G_i : i \in \Bbb N \rangle$ be a sequence
of graphs, each with a Hamilton path.
Let $\Psi ( k , P)$ be the arithmetical sentence
formalizing
``$P$ is a Hamilton path through $G_k$.''
By $\s11$, since
$(\forall k)( \exists P) \Psi (k, P)$, there is a $Y$ such
that $(\forall k) \Psi (k,(Y)_k)$.
Since the desired sequence of paths is
$\Delta_0^1$-definable in $Y$,  it exists by the
recursive comprehension axiom.

The first step in proving that
\therosteritem2 implies \therosteritem1 is to
deduce $\aca$ from \therosteritem2.
Let $f: \Bbb N \rightarrow \Bbb N$ be an injection.
We will show that $Range (f)$ exists.
Construct a sequence of graphs as follows.
Let $v_0 , v_1 , v_2 , \dotsc$ be the vertices
of $G_n$.  Include the edge $(v_0,v_1)$ in $G_n$.
For each $j \in \Bbb N$, if $f(j) \neq n$,
add the edge $(v_{j+1} , v_ {j+2} )$ to $G_n$.
If $f(j) = n$, add $(v_0, v_{j+2} )$ to $G_n$.
By the recursive comprehension axiom,
the sequence
$\langle G_i : i \in \Bbb N \rangle$
exists.
Furthermore, for each $n$, $G_n$ has a Hamiltonian path.
In particular,
if $n \notin Range(f)$ the only Hamiltonian path in $G_n$ is
$v_0 v_1 v_2 \dotso$, while if $f(m)=n$, the only path is
given by $v_{m+1} v_m \dotso v_0 v_{m+2} v_{m+3} \dotso$.
Applying \therosteritem2, we obtain a sequence of paths
$\langle P_i : i \in \Bbb N \rangle$,
and by the recursive comprehension axiom,
the set
$$
Range (f) = \{ n: v_0 ~ \text{is not the first vertex in} ~ P_n \}
$$
exists.  By Theorem 2, this suffices to prove $\aca$.

To complete the deduction of $\s11$ from \therosteritem2,
suppose that $\Psi$ is a $\Sigma_1^1$ formula and
$(\forall k)( \exists X)\Psi (k,X)$.  By Lemma 3.14
of \cite{5}, there is a sequence of trees
$\langle T_i : i \in \Bbb N \rangle$ such that for all
$k \in \Bbb N$,  $P$ is a path through
$T_k$ if and only if $\Psi (k,X)$ , where $X$ is
uniformly $\Delta_1^0$-definable in $P$.
(Lemma 3.14 of \cite{5} has $\aca$ as a hypothesis.)
Let
$\langle G_i : i \in \Bbb N \rangle$
be the sequence of graphs obtained by
applying Lemma 25
to $\langle T_i : i \in \Bbb N \rangle$.
Since $(\forall k)( \exists X )\Psi (k,X)$,
each of the trees has an infinite path, so
each of the graphs has a Hamilton path.
By \therosteritem2 , there is a sequence of paths
$\langle P_i : i \in \Bbb N \rangle$ for the graphs.
By Lemma 25, there is a sequence of paths
$\langle P_i^\prime : i \in \Bbb N \rangle$
for the trees.  Using these paths as a parameter,
arithmetic comprehension suffices to prove the
existence of a set $Y$ such that
$(\forall k )\Psi (k , (Y)_k )$.  Thus
$\s11$ holds, as desired.  \qed
\enddemo

>From Theorem 26, we can draw the following recursion theoretic conclusion.

\proclaim{Corollary 27}
If $\langle G_i : i \in \Bbb N \rangle$
is a hyperarithmetical sequence of graphs, each of which
has a Hamilton path, then there is a hyperarithmetical
sequence $\langle P_i : i \in \Bbb N \rangle$ such that
for each $i$, $P_i$ is a Hamilton path through
$G_i$.
\endproclaim

\demo{Proof}
$\omega$ together with the hyperarithmetical sets is
a model of $\s11$ \cite{16}.  \qed
\enddemo

Using Theorem 24, it is easy to prove:

\proclaim{Theorem 28 ($\rca$)}
The following are equivalent:
\roster
\item
$\atr$. 
\item
If $\langle G_i : i \in  \Bbb N \rangle$ is a sequence of graphs
each of which has at most one Hamilton path, then there is a
set $Z \subseteq \Bbb N$ such that for all $i \in \Bbb N$,
$i \in Z$ if and only if $G_i$ has a Hamilton path.
\endroster
\endproclaim

\demo{Proof}
To prove that \therosteritem1 implies \therosteritem2, apply the
scheme in part \therosteritem2 of Theorem 24,
using ``$X$ is a Hamilton path through $G_i$''
for $\Psi (i,X)$.

To prove the converse, it suffices to deduce part \therosteritem3
of
Theorem 24 using \therosteritem2.
Let $\langle T_i : i \in  \Bbb N \rangle$ be a sequence of trees, each
with at most one path.
Lemma 25 yields a corresponding sequence of graphs, each with at most
one Hamilton path.  The set $Z$ obtained by applying \therosteritem2
satisfies part \therosteritem3 of
Theorem 24.  \qed
\enddemo

The following corollary is a recursion theoretic consequence of Theorem 28.

\proclaim{Corollary 29}
There is a hyperarithmetical sequence of graphs
$\langle G_i : i \in \Bbb N \rangle$, each of
which has at most one Hamilton path, such that
the set
$\{ i \in \Bbb N : G_i~\text{has a Hamilton path} \} $
is not hyperarithmetical.
\endproclaim

\demo{Proof}
$\omega$ together with the hyperarithmetical sets is not a
model of $\atr$ \cite{16}.  \qed
\enddemo

Now we will analyze the third and most difficult task.
Theorem 30 is closely related
to {\smc Harel's} proof \cite{7} that the problem of finding
a Hamiltonian path is $\Sigma^1_1$ complete.

\proclaim{Theorem 30 ($\rca$)}
The following are equivalent:
\roster
\item
$\p11$. 
\item
If $\langle G_i : i \in  \Bbb N \rangle$
is a sequence of graphs, then there is a set $Z \subseteq \Bbb N$
such that $i \in Z$ if and only if
$G_i$ has a Hamilton path.
\endroster
\endproclaim

\demo{Proof}
To prove that \therosteritem1 implies \therosteritem2, assume
\therosteritem1 and let $\langle G_i : i \in  \Bbb N \rangle$
be a sequence of graphs.  By $\Pi_1^1$ comprehension, the set
$$
C = \{ i \in \Bbb N : G_i ~\text{does not have a Hamilton path}\}
$$
exists.  By the recursive comprehension axiom, the desired set $Z$,
which is the complement of $C$, also exists.

To prove the converse, we will use \therosteritem2 to prove that
$\{ n\in \Bbb N : \Psi (n) \}$ exists, where $\Psi (n)$ is a
$\Pi^1_1$ formula.  Note that $\lnot \Psi (n)$ is a $\Sigma_1^1$
formula.  By Lemma 3.14 of \cite{5}, there is a sequence of
trees $\langle T_i : i \in  \Bbb N \rangle$ such that
$T_i$ has a path if and only if $\lnot \Psi (i)$.
By Lemma 25, there is a sequence of associated graphs
$\langle G_i : i \in  \Bbb N \rangle$ such that $G_i$ has
a Hamilton graph if and only if $\lnot \Psi (i)$.
Applying \therosteritem2 yields the set
$Z = \{ n \in \Bbb N : \lnot \Psi (n) \}$.
By the recursive comprehension axiom, the complement of
Z, $\{ n \in \Bbb N : \Psi (n) \}$, also exists.  \qed
\enddemo

Theorem 30 contrasts nicely with Theorem 21.  Since $\p11$ is a
much stronger axiom system than $\aca$, we can conclude that
it is more difficult to determine if certain graphs have Hamilton paths
than to determine if they have Euler paths.  Determining which finite
graphs have Hamilton paths is an NP-complete problem,
while determining which finite graphs have Euler paths is
polynomial time computable.  It would be nice to know if this
sort of parallel is common, and exactly what it signifies.


\enddocument